\documentclass[12pt]{amsart}

\usepackage{amsfonts,amsmath,amsthm}
\usepackage{latexsym}

\theoremstyle{remark}
\newtheorem{remark}{Remark}

\def \a {\alpha}

\def \< {\langle}
\def \> {\rangle}
\def \ra {\rightarrow}

\def \co {4+\frac{\pi^2}{4}}
\def \ct {6+\frac{\pi^2}{6}-\frac{\pi^4}{27\cdot 24\cdot 7}}

\begin{document}

\title[Counterexample to Trotter formula]{Counterexample to the Trotter product
formula for projections}
\author{M\'at\'e~Matolcsi}
\email{matomate@cs.elte.hu}
\address{ELTE TTK, Department of Applied Analysis, 1053 Budapest, Kecskem\'{e}ti
u. 10-12, Hungary.}
\author{Roman~Shvidkoy}
\email{shvidkoy@math.missouri.edu}
\address{Mathematics Department, University of Missouri-Columbia, Columbia, MO
65211, USA.}  \subjclass{47d03}
\thanks{The research of the second author was supported by the
student part of the NSF grant DMS-9800027}
\date{\today}
\maketitle

\begin{abstract}
We constructed a unitary semigroup $(e^{tA})_{t \geq 0}$ on a
Hilbert space and an orthogonal projection $P$ such that the limit
$$
\lim_{n \rightarrow \infty} \left[ e^{\frac{t}{n}A}P \right ]^n
$$
does not exist strongly. A similar example with a positive
contractive semigroup and positive contractive projection on $L_p$
is also constructed.
\end{abstract}

\section{Introduction}
In this short note we construct counterexamples to the following
conjecture.

Suppose $P$ is an orthogonal projection in a Hilbert space $H$ and
$(e^{tA})_{t \geq 0}$ is a contractive semigroup on $H$. Does the
limit
\begin{equation}\label{main}
\lim_{n \rightarrow \infty} \left[ e^{\frac{t}{n}A}P \right ]^n
\end{equation}
exist strongly, for all $t>0$? Or, let $(e^{tA})_{t \geq 0}$ is
positive contractive on $L_p$, for $1<p<\infty$, and $P$ is a
positive contractive projection. Is the same assertion true?

Even though the limit \eqref{main} fails to exist in the above
settings, there are a few important cases when \eqref{main} does
make sense and defines a degenerate semigroup. Such is, for
example, the case when $A$ is associated with a sectorial
sesquilinear form. The reader will find the detailed exposition of
the subject and further references in \cite{Arendt-Ulm97}, where,
in particular, the conjecture was originally posed.

\section{Counterexamples}
\subsection{Hilbert case}
Let us remark that by using the theory of unitary dilations of
contractive semigroups in Hilbert spaces one can reduce the first
question of the conjecture to the case of unitary semigroups.
Therefore, we are looking for a counterexample among unitary
semigroups instead of arbitrary contractive ones.

We carry out our construction in the space $L_2[0,1]$. As an
example of unitary semigroup we take the semigroup of
multiplication by $e^{ith}$, where $h$ is a real-valued,
measurable function on $[0,1]$, to be specified later. We choose
$P$ to be the one-dimensional orthogonal projection onto the space
of constant functions, i.e. $Pf=\mathbf{1}\cdot \int_0^1f(x)dx$.
As a test function on which \eqref{main} will fail for $t=1$, we
take $\mathbf{1}$.

Denoting $c_n=\int_0^1e^{i\frac{1}{n}h(x)}dx$,  the function
$\left[ e^{\frac{1}{n}A}P \right ]^n (\mathbf{1})$ becomes
$c_n^{n-1}e^{i\frac{1}{n}h}$. However, by the Lebesgue Dominated
Convergence Theorem, $\lim_{n \rightarrow \infty} c_n =1$ as well
as $\lim_{n \rightarrow \infty}e^{i\frac{1}{n}h}=\mathbf{1}$ in
$L_2[0,1]$. So, $\lim_{n \rightarrow \infty}\left[
e^{\frac{1}{n}A}P \right ]^n (\mathbf{1})$ exists in $L_2[0,1]$ if
and only if the numerical limit
\begin{equation}\label{limit}
\lim_{n \rightarrow \infty} c_n^n
\end{equation}
exists. Now we specify the function $h$, for which we prove that
\eqref{limit} diverges. Put
$h=\sum_{k=1}^{\infty}\chi_{(1/2^k,1/2^{k-1}]}2^k \pi$. Then
$c_n=\sum_{k=1}^{\infty}\frac{1}{2^k}e^{i\frac{1}{n}2^k \pi}$. We
show the following two inequalities
\begin{eqnarray}
\label{first} \liminf_{n \rightarrow \infty} |c_{2^n}|^{2^n} \geq
e^{-(\co)} \\
\label{second} \limsup_{n \rightarrow \infty} |c_{2^n 3}|^{2^n 3}
\leq e^{-(\ct)}.
\end{eqnarray}
Noticing that $\co<\ct$ we get the desired result.

Let us show \eqref{first} first. Observe that
$$
c_{2^n}=\sum_{k=1}^{n-1}\frac{1}{2^k}e^{i\frac{2^k}{2^n}
\pi}-\frac{1}{2^n}+\sum_{k=n+1}^{\infty}\frac{1}{2^k}=\sum_{k=1}^{n-1}\frac{1}{2^k}e^
{i\frac{2^k}{2^n} \pi}.
$$
Using the inequality $\cos(\a) \geq 1- \frac{\a^2}{2}$ we get
\begin{eqnarray*}
|c_{2^n}|\geq
|\Re{c_{2^n}}|=\sum_{k=1}^{n-2}\frac{1}{2^k}\cos(\frac{2^k}{2^n}\pi)
\geq
\sum_{k=1}^{n-2}\frac{1}{2^k}(1-\frac{\pi^2}{2}\frac{4^k}{4^n})
\\= 1-\frac{4}{2^n}-\frac{\pi^2}{2}\frac{1}{4^n}(2^{n-1}-2)=
1-\frac{1}{2^n}(\co)+\frac{\pi^2}{4^n}.
\end{eqnarray*}
Since $\lim_{N \ra \infty}(1+\frac{a}{N}+\frac{b}{N^2})^{N}=e^a$,
we obtain \eqref{first}.

To prove \eqref{second} let us simplify $c_{2^n3}$. We have
\begin{eqnarray*}
c_{2^n3}=\sum_{k=1}^{n-1}\frac{1}{2^k}e^{i\frac{2^k}{2^n3}
\pi}+\frac{1}{2^n}e^{i\frac{1}{3}\pi}+\sum_{k=n+1}^{\infty}\frac{1}{2^k}e^{i\frac{2^{k-n}}{3}\pi}\\
=\sum_{k=1}^{n-1}\frac{1}{2^k}e^{i\frac{2^k}{2^n3}
\pi}+\frac{1}{2^n}(\frac{1}{2}+i\frac{\sqrt{3}}{2})+\frac{1}{2^n}
\sum_{k=1}^{\infty}\frac{1}{2^k}e^{i\frac{2^{k}}{3}\pi}.
\end{eqnarray*}
Notice that
$e^{i\frac{2^{k}}{3}\pi}=e^{i(-1)^{k+1}\frac{2}{3}\pi}=-\frac{1}{2}+i(-1)^{k+1}\frac{\sqrt{3}}{2}$.
Thus,
$\sum_{k=1}^{\infty}\frac{1}{2^k}e^{i\frac{2^{k}}{3}\pi}=-\frac{1}{2}+i\frac{\sqrt{3}}{6}$.
After these computations $c_{2^n3}$ becomes
$$
\sum_{k=1}^{n-1}\frac{1}{2^k}e^{i\frac{2^k}{2^n3} \pi} +
i\frac{2\sqrt{3}}{2^n3}.
$$
Now using the inequality $\cos(\a)\leq
1-\frac{\a^2}{2}+\frac{\a^4}{24}$ we obtain the following estimate
\begin{eqnarray*}
|\Re{c_{2^n3}}| \leq \sum_{k=1}^{n-1}\frac{1}{2^k}\left(
1-\frac{\pi^2}{18}\frac{4^k}{4^n}+\frac{\pi^4}{81\cdot
24}\frac{16^k}{16^n} \right) \\
= 1 -
\frac{1}{2^{n-1}}-\frac{\pi^2}{18}\frac{2^n-2}{4^n}+\frac{\pi^4}{81\cdot
24}\frac{8^n-8}{16^n 7}\\
= 1- \frac{1}{2^n3}\left( \ct \right)+ \frac{a}{(2^n3)^2} +
\frac{b}{(2^n3)^4},
\end{eqnarray*}
for some constants $a$ and $b$. Similarly, using $\sin(\a)\leq
\a$, we have
$$
|\Im{c_{2^n3}}| \leq
\sum_{k=1}^{n-1}\frac{1}{2^k}\frac{2^k}{2^n3}\pi+\frac{2\sqrt{3}}{2^n3}\leq
\frac{(n+1)\pi}{2^n3}.
$$
Thus,
\begin{eqnarray*}
|c_{2^n3}|^{2^n3}&=&
\left(|\Re{c_{2^n3}}|^2+|\Im{c_{2^n3}}|^2\right)^{\frac{2^n3}{2}} \\
&\leq&  ( 1- \frac{2}{2^n3}(\ct) + (\frac{2}{2^n3})^2(n+1)^2 a_1 \\
&+& (\frac{2}{2^n3})^2 a_2 + \ldots +(\frac{2}{2^n3})^8 a_8
)^{\frac{2^n3}{2}}.
\end{eqnarray*}
Passing to the upper limit as $n \ra \infty$, we finally obtain
\eqref{second}.

\subsection{$L_p$-case}

Our second example is on the Hilbert space $L_2[0,2\pi]$, but now
for a positive contractive semigroup and positive contractive
projection.

We take $e^{tA}f(x)=f(x+2\pi t)$, regarding $f$ as a
$2\pi$-periodic function. Now let $P$ be the orthogonal projection
onto the space spanned by the positive norm-one function
$g(x)=\frac{1}{\sqrt{34\pi}} \left[ 4+ \sum_{k=0}^{\infty}
\frac{1}{\sqrt{2^k}}\cos{2^k x}\right]$. Notice that like in the
previous example our projection is one-dimensional (see Remark
\ref{rem} below). Simple substitution shows that \eqref{main}
evaluated at $g$ for $t=1$ exists if and only if the numerical
limit
$$
\lim_{n\ra \infty}
\left[\int_0^{2\pi}g(x)g(x+\frac{1}{n})dx\right]^{n}
$$
exists. Denoting
$$
c_n= \int_0^{2\pi}g(x)g(x+\frac{1}{n})dx
$$
and using the orthogonality of cosines, we get
$$
c_n=\frac{16}{17}+\frac{1}{17}\sum_{k=1}^{\infty}\frac{1}{2^k}\cos{\frac{2^k}{n}\pi}
$$
Now following the same calculations as for the first example, we
obtain inequalities \eqref{first} and \eqref{second} with powers
doubled on the right hand sides.

This disproves the second conjecture.

\begin{remark}\label{rem} As we have already noticed the
projections in our examples are one-dimensional. It would be
interesting to know what property of a semigroup on a Hilbert
space is responsible for the existence of \eqref{main} for all
one-dimensional, or more specifically, one-dimensional orthogonal
projections.
\end{remark}

The authors are grateful to Wolfgang~Arendt, Andr\'as~B\'atkai and
B\'alint~Farkas for helpful conversations.

\end{document}